\documentclass{article}                                     
                                                                             
\input{amssym.def}                                                           
\input{amssym.tex}                                                           
                                                                             
\begin{document}                                                             
\title{On a Teichm\"uller functor between the categories of 
complex tori and the Effros-Shen algebras}

\author{Igor V. Nikolaev
\footnote{Partially supported 
by NSERC.}}


\date{}
 \maketitle


\newtheorem{thm}{Theorem}
\newtheorem{lem}{Lemma}
\newtheorem{dfn}{Definition}
\newtheorem{rmk}{Remark}
\newtheorem{cor}{Corollary}
\newtheorem{cnj}{Conjecture}
\newtheorem{exm}{Example}


\begin{abstract}
A covariant functor from the category of the complex tori to the category
of the Effros-Shen algebras is constructed. The functor maps isomorphic
complex tori to the stably isomorphic  Effros-Shen algebras.
Our construction is based on the Teichm\"uller theory of
the Riemann surfaces.

\vspace{7mm}


\vspace{5mm}
{\it AMS (MOS) Subj. Class.:  14H52, 46L85}
\end{abstract}

\noindent
{\bf A. The complex tori.} 
Let $\omega_1$ and $\omega_2$ be a pair of the non-zero complex numbers,  
which are linearly independent over ${\Bbb R}$.  Consider a lattice 
$\Lambda:={\Bbb Z}\omega_1+{\Bbb Z}\omega_2$ in the complex plane ${\Bbb C}$ and
the quotient space  ${\Bbb C}/\Lambda$. The space ${\Bbb C}/\Lambda$ 
is known as a {\it complex torus}.  It is easy to see, that 
the conformal transformation of the complex plane $z\to\pm {\omega_2\over\omega_1}z$
brings the complex torus  to a normal form ${\Bbb C}/\Lambda_{\tau}$,
where $\Lambda_{\tau}:={\Bbb Z}+{\Bbb Z}\tau$ and  $\tau\in {\Bbb H}:=\{z\in {\Bbb C}~|~Im (z)>0\}$.
The transformation law in a lattice implies, that the complex tori
 ${\Bbb C}/\Lambda_{\tau}$ and   ${\Bbb C}/\Lambda_{\tau'}$
are conformally equivalent (isomorphic),  whenever $\tau'\equiv \tau ~mod~SL_2({\Bbb Z})$, 
i.e. $\tau'= {a +b\tau\over c+d\tau}$,   where $a,b,c,d\in {\Bbb Z}$ and $ad-bc=1$.

\medskip\noindent
{\bf B. The Effros-Shen algebras.}
 Let $\theta>0$ be an irrational number given by the regular
continued fraction $[a_0, a_1, a_2,\dots]$, where $a_0\in {\Bbb N}\cup 0$
and $a_i\in {\Bbb N}$ for $i\ge 1$.  
By an {\it Effros-Shen algebra} \cite{EfSh1},  one understands
the $AF$-algebra ${\Bbb A}_{\theta}$ given by the Bratteli diagram:

\begin{figure}[here]
\begin{picture}(300,60)(0,0)
\put(110,30){\circle{3}}
\put(120,20){\circle{3}}
\put(140,20){\circle{3}}
\put(160,20){\circle{3}}
\put(120,40){\circle{3}}
\put(140,40){\circle{3}}
\put(160,40){\circle{3}}

\put(110,30){\line(1,1){10}}
\put(110,30){\line(1,-1){10}}
\put(120,42){\line(1,0){20}}
\put(120,40){\line(1,0){20}}
\put(120,38){\line(1,0){20}}
\put(120,40){\line(1,-1){20}}
\put(120,20){\line(1,1){20}}
\put(140,41){\line(1,0){20}}
\put(140,39){\line(1,0){20}}
\put(140,40){\line(1,-1){20}}
\put(140,20){\line(1,1){20}}

\put(180,20){$\dots$}
\put(180,40){$\dots$}

\put(125,52){$a_0$}
\put(145,52){$a_1$}

\end{picture}

\caption{The Effros-Shen algebra ${\Bbb A}_{\theta}$.}
\end{figure}
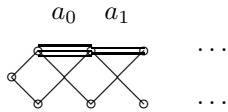

\noindent
where $a_i$ indicate the number of edges in the upper row of the graph.
Recall that the Effros-Shen algebras ${\Bbb A}_{\theta}, {\Bbb A}_{\theta'}$ are 
said to be  stably isomorphic (Morita equivalent), if  ${\Bbb A}_{\theta}\otimes {\cal K}
\cong {\Bbb A}_{\theta'}\otimes {\cal  K}$,  where ${\cal K}$ is a $C^*$-algebra
of the compact operators. A remarkable result, proved in \cite{EfSh1}, says that  
${\Bbb A}_{\theta}$ and  ${\Bbb A}_{\theta'}$ are stably isomorphic  
if and only if $\theta'\equiv \theta~mod~SL_2({\Bbb Z})$.

\medskip\noindent
{\bf C. Motivation and background.}  
Comparing the category of the complex tori with such of the Effros-Shen 
algebras, one cannot fail to observe that for the generic objects in the respective
categories, the corresponding morphisms (modulo the inner automorphisms) are isomorphic as groups.  
Assuming that  the observation is not a simple coincidence,  one can ask the following 
question.

\medskip\noindent
{\sf Main problem.} 
{\it Let ${\cal A}$ be the category of the complex tori, described in item (A), and let ${\cal B}$
be the category of the Effros-Shen algebras, described in item (B).  Construct a functor 
(if any) $F: {\cal A}\to {\cal B}$, which maps isomorphic complex tori to the stably isomorphic
Effros-Shen algebras.}

\medskip\noindent
The question attracted attention of  both the algebraic geometers and the operator
algebraists. Manin \cite{Man1} and Soibelman \cite{Soi1}, \cite{Soi2} were apparently
the first to study such a functor.  (Note that the authors usually consider a category of the 
noncommutative tori \cite{Rie1},  i.e. the universal $C^*$-algebras generated by the two unitary operators $U$ and $V$,
which satisfy the commutation relation $VU=e^{2\pi i\theta}UV$.  However, it is well known that  
the two objects are closely related \cite{PiVo1}.)  The topic was pursued by 
Polishchuk \cite{Pol1} -- \cite{Pol4} and Polishchuk-Schwarz \cite{PoSch1},  using the methods of the
 homological algebra and the algebraic geometry. The works of Kontsevich \cite{Kon1}
and Soibelman-Vologodsky \cite{SoVo1} develop the ideas of a homological  mirror symmetry
and the deformation quantization of the elliptic curves.  Finally, Mahanta \cite{Mah1},
Mahanta-van Suijlekom \cite{MaSui1}, Plazas \cite{Pla1}, \cite{Pla2}
and Taylor \cite{Tay1}, \cite{Tay2} elaborated the ideas of Polishchuk-Schwarz
and Manin, respectively.

\medskip\noindent
{\bf D. The measured foliations.} 
 A measured foliation, ${\cal F}$, on a surface $X$
is a  partition of $X$ into the singular points $x_1,\dots,x_n$ of
the order $k_1,\dots, k_n$ and the regular leaves (the $1$-dimensional submanifolds). 
On each  open cover $U_i$ of $X-\{x_1,\dots,x_n\}$ there exists a non-vanishing
real-valued closed 1-form $\phi_i$  such that: 
(i)  $\phi_i=\pm \phi_j$ on $U_i\cap U_j$;
(ii) at each $x_i$ there exists a local chart $(u,v):V\to {\Bbb R}^2$
such that for $z=u+iv$, it holds $\phi_i=Im~(z^{k_i\over 2}dz)$ on
$V\cap U_i$ for some branch of $z^{k_i\over 2}$. 
The pair $(U_i,\phi_i)$ is called an atlas for the measured foliation ${\cal F}$.
Finally, a measure $\mu$ is assigned to each segment $(t_0,t)\in U_i$, which is  transverse to
the leaves of ${\cal F}$, via the integral $\mu(t_0,t)=\int_{t_0}^t\phi_i$. The 
measure is invariant along the leaves of ${\cal F}$, hence the name. 
Note that in the case $X=T^2$ is a two-dimensional torus (our main concern),
every measured foliation is given by  a family of the parallel lines of a slope $\theta>0$,
see Fig. 2.

\begin{figure}[here]
\begin{picture}(300,60)(-30,0)

\put(130,10){\line(1,0){40}}
\put(130,10){\line(0,1){40}}
\put(130,50){\line(1,0){40}}
\put(170,10){\line(0,1){40}}

\put(130,40){\line(2,1){20}}
\put(130,30){\line(2,1){40}}
\put(130,20){\line(2,1){40}}
\put(130,10){\line(2,1){40}}

\put(150,10){\line(2,1){20}}

\end{picture}

\caption{A measured foliation on the torus ${\Bbb R}^2/{\Bbb Z}^2$.}
\end{figure}
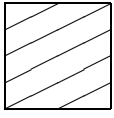

\medskip\noindent
{\bf E.  The Hubbard-Masur homeomorphism.} 
Let $T(g)$ be the Teichm\"uller space of the topological surface $X$ of genus $g\ge 1$,
i.e. the space of the complex structures on $X$. 
Consider the vector bundle $p: Q\to T(g)$ over $T(g)$,  whose fiber above a point 
$S\in T_g$ is the vector space $H^0(S,\Omega^{\otimes 2})$.   
Given a non-zero $q\in Q$ above $S$, we can consider the horizontal measured foliation
${\cal F}_q\in \Phi_X$ of $q$, where $\Phi_X$ denotes the space of the equivalence
classes of the measured foliations on $X$. If $\{0\}$ is the zero section of $Q$,
the above construction defines a map $Q-\{0\}\longrightarrow \Phi_X$. 
For any ${\cal F}\in\Phi_X$, let $E_{\cal F}\subset Q-\{0\}$ be the fiber
above ${\cal F}$. In other words, $E_{\cal F}$ is a subspace of the holomorphic 
quadratic forms,  whose horizontal trajectory structure coincides with the 
measured foliation ${\cal F}$. 
Note that, if ${\cal F}$ is a measured foliation with the simple zeroes (a generic case),  
then $E_{\cal F}\cong {\Bbb R}^n - 0$, while $T(g)\cong {\Bbb R}^n$, where $n=6g-6$ if
$g\ge 2$ and $n=2$ if $g=1$.

\medskip\noindent
{\bf Theorem (Hubbard-Masur \cite{HuMa1})}
{\it The restriction of $p$ to $E_{\cal F}$ defines a homeomorphism (an embedding)
$h_{\cal F}: E_{\cal F}\to T(g)$.}

\medskip\noindent
{\bf F.  The Teichm\"uller space and measured foliations.}
The Hubbard-Masur result implies that the measured foliations  parametrize  
the space $T(g)-\{pt\}$, where $pt= h_{\cal F}(0)$.
Indeed, denote by  ${\cal F}'$ a vertical trajectory structure of  $q$. Since ${\cal F}$
and ${\cal F}'$ define $q$, and ${\cal F}=Const$ for all $q\in E_{\cal F}$, one gets a homeomorphism 
between $T(g)-\{pt\}$ and  $\Phi_X$, where $\Phi_X\cong {\Bbb R}^n-0$ is the space of 
equivalence classes of the measured foliations ${\cal F}'$ on $X$. 
Note that the above parametrization depends on a foliation ${\cal F}$.
However, there exists a unique canonical homeomorphism $h=h_{\cal F}$
as follows. Let $Sp ~(S)$ be the length spectrum of the Riemann surface
$S$ and $Sp ~({\cal F}')$ be the set positive reals $\inf \mu(\gamma_i)$,
where $\gamma_i$ runs over all simple closed curves, which are transverse 
to the foliation ${\cal F}'$. A canonical homeomorphism 
$h=h_{\cal F}: \Phi_X\to T(g)-\{pt\}$ is defined by the formula
$Sp ~({\cal F}')= Sp ~(h_{\cal F}({\cal F}'))$ for $\forall {\cal F}'\in\Phi_X$. 
Thus, the following corollary is true.

\medskip\noindent
{\bf Corollary}
{\it ~There exists a canonical  homeomorphism $h:\Phi_X\to T(g)-\{pt\}$.}

\medskip\noindent
{\bf G.  A parametrization of ${\Bbb H}-\{pt\}$  by the measured foliations.}
In the case $X=T^2$, the picture simplifies.  First,  notice that $T(1)\cong {\Bbb H}$.
Since $q\ne 0$ there are no singular points and each $q\in H^0(S, \Omega^{\otimes 2})$
has the form $q=\omega^2$, where $\omega$ is a nowhere zero  holomorphic differential
on the complex torus $S$. 
(Note that $\omega$ is just a constant times $dz$, and hence its vertical
trajectory structure is just a family of the parallel lines of a slope $\theta$,
see e.g. Strebel \cite{S}, pp. 54--55.) 
Therefore,  $\Phi_{T^2}$ consists of the equivalence classes of the non-singular 
measured foliations on the two-dimensional torus.  It is well known (the Denjoy theory), 
that every such foliation  is measure equivalent to the foliation of a slope $\theta$ and a 
transverse measure $\mu>0$, which is  invariant along the leaves of the foliation (Fig.2). 
Thus,  one obtains  a canonical bijection $h: \Phi_{T^2}\to {\Bbb H}-\{pt\}$.

\medskip\noindent
{\bf H.  The  lattices.} 
Let ${\Bbb C}$ be the complex plane. 
A {\it lattice} is a triple $(\Lambda, {\Bbb C}, j)$,  where
$\Lambda\cong {\Bbb Z}^2$ and $j: \Lambda\to {\Bbb C}$ is an injective 
homomorphism with the discrete image.  A {\it morphism} of lattices 
$(\Lambda, {\Bbb C}, j)\to (\Lambda', {\Bbb C}, j')$  
is the identity $j\circ\psi=\varphi\circ j'$
where $\varphi$ is a group homomorphism and $\psi$ is a ${\Bbb C}$-linear
map. It is not hard to see, that any isomorphism class of a lattice contains
a representative given by $j: {\Bbb Z}^2\to  {\Bbb C}$ such that $j(1,0)=1,
j(0,1)=\tau\in {\Bbb H}$. The category of lattices, ${\cal L}$, consists of $Ob~({\cal L})$,
which are lattices $(\Lambda, {\Bbb C}, j)$ and morphisms $H(L,L')$
between $L,L'\in Ob~({\cal L})$ which coincide with the morphisms
of lattices specified above. For any  $L,L',L''\in Ob~({\cal L})$
and any morphisms $\varphi': L\to L'$, $\varphi'': L'\to L''$ a 
morphism $\phi: L\to L''$ is the {\it composite} of $\varphi'$ and
$\varphi''$, which we write as $\phi=\varphi''\varphi'$. 
The {\it identity} morphism, $1_L$, is a morphism $H(L,L)$.
Note that the lattices are bijective with the complex tori via the formula
$(\Lambda, {\Bbb C}, j)\mapsto {\Bbb C}/j(\Lambda)$. 
Therefore,  ${\cal L}\cong {\cal A}$.

\medskip\noindent
{\bf I.  The pseudo-lattices.}  
Let ${\Bbb R}$ be the real line. 
A {\it pseudo-lattice} (of rank 2) is a triple $(\Lambda, {\Bbb R}, j)$, where
$\Lambda\cong {\Bbb Z}^2$ and $j: \Lambda\to {\Bbb R}$ is a homomorphism.  
A morphism of the pseudo-lattices $(\Lambda, {\Bbb R}, j)\to (\Lambda', {\Bbb R}, j')$
is  the identity $j\circ\psi=\varphi\circ j'$,
where $\varphi$ is a group homomorphism and $\psi$ is an inclusion 
map (i.e. $j'(\Lambda')\subseteq j(\Lambda)$).  
Any isomorphism class of a pseudo-lattice contains
a representative given by $j: {\Bbb Z}^2\to  {\Bbb R}$, such that $j(1,0)=\lambda_1,
j(0,1)=\lambda_2$, where $\lambda_1,\lambda_2$ are the positive reals.
The pseudo-lattices make up a category, which we denote by ${\cal PL}$.
\begin{lem}\label{lm1}
The pseudo-lattices are bijective with the measured foliations on the torus
via the formula $(\Lambda, {\Bbb R}, j)\mapsto {\cal F}_{\lambda_2/\lambda_1}^{\lambda_1}$,
where  ${\cal F}_{\lambda_2/\lambda_1}^{\lambda_1}$ is a foliation of the slope $\theta=\lambda_2/\lambda_1$
and measure $\mu=\lambda_1$.  
\end{lem}
{\it Proof}. Define a pairing by the formula
$(\gamma, Re~\omega)\mapsto \int_{\gamma} Re~\omega$, where  $\gamma\in H_1(T^2, {\Bbb Z})$
and $\omega\in H^0(S; \Omega)$. The trajectories of the closed differential $\phi:=Re~\omega$ 
define a measured foliation on $T^2$. Thus, in view of the pairing, the linear spaces $\Phi_{T^2}$
and  $Hom~(H_1(T^2, {\Bbb Z}); {\Bbb R})$ are isomorphic. Notice that 
the latter space coincides with  the  space of the pseudo-lattices.
To obtain an explicit bijection formula, let us evaluate the integral:
\begin{equation}\label{eq1}
\int_{{\Bbb Z}\gamma_1+{\Bbb Z}\gamma_2}\phi ={\Bbb Z}\int_{\gamma_1}\phi + 
{\Bbb Z}\int_{\gamma_2}\phi= {\Bbb Z}\int_0^1\mu dx + {\Bbb Z}\int_0^1\mu dy,
\end{equation}
where $\{\gamma_1,\gamma_2\}$ is a basis in $H_1(T^2, {\Bbb Z})$. 
Since ${dy\over dx}=\theta$, one gets:
\begin{equation}\label{eq2}
\left\{
\begin{array}{cccc}
\int_0^1\mu dx  &= \mu  &= \lambda_1  & \nonumber\\
\int_0^1\mu dy  &= \int_0^1\mu\theta dx   &= \mu\theta  &=
\lambda_2. 
\end{array}
\right.
\end{equation}
Thus, $\mu=\lambda_1$ and $\theta={\lambda_2\over\lambda_1}$.
$\square$

\medskip
It follows from lemma \ref{lm1} and the canonical bijection 
 $h: \Phi_{T^2}\to {\Bbb H}-\{pt\}$, that ${\cal L}\cong {\cal PL}$
are the equivalent categories.

\medskip\noindent
{\bf J. The projective pseudo-lattices.}
Finally,  a  {\it projective pseudo-lattice} (of rank 2) is a triple 
$(\Lambda, {\Bbb R}, j)$, where $\Lambda\cong {\Bbb Z}^2$ and $j: \Lambda\to {\Bbb R}$ is a
homomorphism. A morphism of the projective pseudo-lattices
$(\Lambda, {\Bbb C}, j)\to (\Lambda', {\Bbb R}, j')$
is  the identity $j\circ\psi=\varphi\circ j'$,
where $\varphi$ is a group homomorphism and $\psi$ is an ${\Bbb R}$-linear
map. (Notice, that unlike the case of the pseudo-lattices, $\psi$ is a scaling
map as opposite to an inclusion map. Thus,  the two pseudo-lattices
can  be projectively equivalent, while being distinct in the category ${\cal PL}$.) 
It is not hard to see that any isomorphism class of a projective pseudo-lattice 
contains a representative given by $j: {\Bbb Z}^2\to  {\Bbb R}$ such that $j(1,0)=1,
j(0,1)=\theta$, where $\theta$ is a positive real.
Note that the projective pseudo-lattices are bijective with the Effros-Shen algebras,
via the formula $(\Lambda, {\Bbb R}, j)\mapsto {\Bbb A}_{\theta}$.
The projective pseudo-lattices make up a category, which we shall 
denote by ${\cal PPL}$.
\begin{lem}\label{lm2}
${\cal PPL}\cong {\cal B}$.
\end{lem}
{\it Proof.}
An isomorphism $\varphi: \Lambda\to\Lambda'$ acts by the formula
$1\mapsto a+b\theta$, $\theta\mapsto c+d\theta$, where $ad-bc=1$
and $a,b,c,d\in {\Bbb Z}$. Therefore,  $\theta'={c+d\theta\over a+b\theta}=\theta~mod~SL_2({\Bbb Z})$.
Thus, the isomorphic projective pseudo-lattices map to the stably isomorphic Effros-Shen algebras.
$\square$

\medskip\noindent
{\bf K. The map $F$ and main results.} 
To finish the construction of a map $F: {\cal A}\to {\cal B}$, consider
a composition of the following morphisms:
\begin{equation}\label{eq3}
{\cal A}
\buildrel\rm\sim
\over\longrightarrow
{\cal L}
\buildrel\rm\sim
\over\longrightarrow
{\cal PL}
\buildrel\rm F
\over\longrightarrow
{\cal PPL}
\buildrel\rm\sim
\over\longrightarrow
 {\cal B},
\end{equation}
where all the arrows, but $F$, have been defined. 
To define $F$,  let $PL\in {\cal PL}$  be a pseudo-lattice,
such that $PL=PL(\lambda_1,\lambda_2)$,  where $\lambda_1=j(1,0), \lambda_2=j(0,1)$
are positive reals.  Let $PPL\in {\cal PPL}$  be a projective pseudo-lattice,
such that $PPL=PL(\theta)$,  where $j(1,0)=1$ and  $j(0,1)=\theta$ is a positive real.
Then $F: {\cal PL}\to {\cal PPL}$ is given by the formula 
$PL(\lambda_1, \lambda_2)\longmapsto PPL\left({\lambda_2\over\lambda_1}\right)$.  
It is easy to see, that $Ker ~F\cong (0,\infty)$ and $F$ is not an injective map. 
Since all the arrows, but $F$,  in the formula (\ref{eq3}) are the isomorphisms between 
the categories, one gets a map $F: {\cal A}\to {\cal B}$. 
\begin{thm}\label{thm1}
The map $F: {\cal A}\to {\cal B}$  is a  covariant non-injective functor with $Ker ~F\cong (0,\infty)$, 
which maps isomorphic complex tori to the stably isomorphic Effros-Shen algebras.
\end{thm}
{\it Proof.} (i) Let us show that $F$ maps isomorphic
complex tori to the stably isomorphic Effros-Shen algebras. 
Let  ${\Bbb C}/({\Bbb Z}\omega_1 +{\Bbb Z}\omega_2)$
be a complex torus.  Recall that the periods 
$\omega_1=\int_{\gamma_1}\omega_E$ and $\omega_2=\int_{\gamma_2}\omega_E$,
where $\omega_E=dz$ is an invariant (N\'eron) differential on the complex
torus and $\{\gamma_1,\gamma_2\}$ is a basis in $H_1(T^2, {\Bbb Z})$.
The map $F$ can be written as: 
\begin{equation}\label{eq4}
{\Bbb C}/\Lambda_{(\int_{\gamma_2}\omega_E) / (\int_{\gamma_1}\omega_E)}
\buildrel\rm F\over
\longmapsto
{\Bbb A}_{(\int_{\gamma_2}\phi)/(\int_{\gamma_1}\phi)},
\end{equation}
where $\phi=Re ~\omega$ is a closed differential defined earlier.
Note that every isomorphism in the category ${\cal A}$ is induced by 
an orientation preserving automorphism, $\varphi$, of the torus $T^2$. 
The action of $\varphi$ on the homology basis $\{\gamma_1,\gamma_2\}$
of $T^2$ is given by the formula:
\begin{equation}\label{eq5}
\left\{
\begin{array}{cc}
\gamma_1' &= a\gamma_1+b\gamma_2\nonumber\\
\gamma_2' &= c\gamma_1+d\gamma_2
\end{array}
\right.
\hbox{, ~~ where}
\quad\left(\matrix{a & b\cr c & d}\small\right)\in SL_2({\Bbb Z}). 
\end{equation}
The functor $F$ acts   by the formula:
\begin{equation}\label{eq6}
\tau={\int_{\gamma_2}\omega_E\over\int_{\gamma_1}\omega_E}
\longmapsto
\theta={\int_{\gamma_2}\phi\over\int_{\gamma_1}\phi}.
\end{equation}

\smallskip\noindent
(a) From the left-hand side of (\ref{eq6}), one obtains 
\begin{equation}\label{eq7}
\left\{
\begin{array}{ccccc}
\omega_1' &= \int_{\gamma_1'}\omega_E &=  \int_{a\gamma_1+b\gamma_2}\omega_E  &= 
a\int_{\gamma_1}\omega_E +b\int_{\gamma_2}\omega_E &= a\omega_1+b\omega_2\nonumber\\
\omega_2' &= \int_{\gamma_2'}\omega_E &=  \int_{c\gamma_1+d\gamma_2}\omega_E  &= 
c\int_{\gamma_1}\omega_E +d\int_{\gamma_2}\omega_E &= c\omega_1+d\omega_2,
\end{array}
\right.
\end{equation}
and therefore $\tau'={\int_{\gamma_2'}\omega_E\over\int_{\gamma_1'}\omega_E}=
{c+d\tau\over a+b\tau}$.

\smallskip\noindent
(b) From the right-hand side of (\ref{eq6}), one obtains
\begin{equation}\label{eq8}
\left\{
\begin{array}{ccccc}
\lambda_1' &= \int_{\gamma_1'}\phi &=  \int_{a\gamma_1+b\gamma_2}\phi  &= 
a\int_{\gamma_1}\phi +b\int_{\gamma_2}\phi &= a\lambda_1+b\lambda_2\nonumber\\
\lambda_2' &= \int_{\gamma_2'}\phi &=  \int_{c\gamma_1+d\gamma_2}\phi  &= 
c\int_{\gamma_1}\phi +d\int_{\gamma_2}\phi &= c\lambda_1+d\lambda_2,
\end{array}
\right.
\end{equation}
and therefore $\theta'={\int_{\gamma_2'}\phi\over\int_{\gamma_1'}\phi}=
{c+d\theta\over a+b\theta}$.  Comparing (a) and (b),  one concludes
that $F$ maps  the isomorphic complex tori to the stably isomorphic 
Effros-Shen algebras.

\medskip
(ii) Let us show that $F$ is a covariant functor, i.e. $F$ does not reverse 
the arrows. Indeed, it can be verified directly using the above  formulas, that 
$F(\varphi_1\varphi_2)=\varphi_1\varphi_2=F(\varphi_1)F(\varphi_2)$
for any pair of the isomorphisms $\varphi_1,\varphi_2\in Aut~(T^2)$. 
Theorem \ref{thm1} is proved.
$\square$

\bigskip\noindent
{\sf Acknowledgments.} 
I wish to thank G. A. Elliott, Yu. I. Manin and L. D. Taylor
for helpful discussions.   The referee suggestions are kindly acknowledged
and incorporated in the text.



\vskip1cm

\textsc{The Fields Institute for Mathematical Sciences, Toronto, ON, Canada,  
E-mail:} {\sf igor.v.nikolaev@gmail.com}

\smallskip
{\it Current address: 101-315 Holmwood Ave., Ottawa, ON, Canada, K1S 2R2}

\end{document}